\documentclass[11pt,twoside]{article}

\date{}
\usepackage{epsfig}
\begin{document}
\begin{figure}[h]
{\Large{\bf  {\emph{Mathematical Sciences}}}}
\hspace{1.1cm}\hspace{0.9cm}%\includegraphics[scale=0.15]{Logo.png}
\hspace{1.1cm}{\small{\bf  {\emph{Vol. 1, No. 1,2 (2007) 01-12}}}}
\end{figure}
\begin{center}
{\Large \bf   Differential Invariants of SL(2) and SL(3)-actions
\\ on ${\bf   R}^2$ Differential Invariants of \\ SL(2) and
SL(3)-actions on ${\bf   R}^2$}

{\small \bf   Mehdi Nadjafikhah$^{a,}$\footnote{\footnotesize
Corresponding Author. E-mail Address: m$\_$nadjafikhah@iust.ac.ir
}, Seyed Reza Hejazi$^b$}
\end{center}
{\tiny \bf   $^a$Fact. of Math., Dept. of Pure Math., Iran University of Science and Technology, Narmak, Tehran, I.R. Iran.}\\
{\tiny \bf   $^b$Same address.}
\newtheorem{Theorem}{\quad Theorem}[section]
\newtheorem{Definition}[Theorem]{\quad Definition}
\newtheorem{Corollary}[Theorem]{\quad Corollary}
\newtheorem{Lemma}[Theorem]{\quad Lemma}
\newtheorem{Example}[Theorem]{\quad Example}
\begin{abstract}
The main purpose of this paper is calculation of differential
invariants which arise from prolonged actions of two Lie groups
SL(2) and SL(3) on the nth jet space of ${\bf    R}^2$. It is
necessary to calculate nth prolonged infinitesimal generators of
the action.
\date{Dec. 18, 2007}
\\
\copyright {\small \hspace{0.15cm}200x Published by Islamic AZAD
University-Karaj Branch.}
\end{abstract}
\section*{Introduction}
Differential invariants theory is one of the most important
concept in differential equations theory and differential
geometry. The study of this theory will help us to analyze
applications of geometry in differential equations. In this paper
we will study some properties of two Lie group actions SL(2) and
SL(3). After finding the Lie algebras of each groups, we will
calculate all infinitesimal generators and their nth
prolongations. Next we are going to find the number of
functionally independent differential invariants up to order n
$(i_n)$, the number of strictly functionally independent
differential invariants of order n $(j_n)$, generic orbit
dimensions of nth prolonged action $(s_n)$ and isotropic group
dimensions of prolonged action $(h_n)$. It is necessary to realize
that the symbol $u_{(k)}$ means kth derivative of $u$ with respect
to $x$.
\section{Jet and Prolongation}
This section is devoted to study of proper geometric context for
"jet spaces" or "jet bundles", well known to 19th century
practitioners, but first formally defined by Ehresmann.
Prolongation will defined after study of jet to find the
differential invariants.
\subsection{Introduction}
In this part we are going to define some fundamental concepts of
differential equations theory. The two most important are Jet and
Prolongation but it needs to analyze some definitions of
differential equations.
\begin{Definition}
Suppose $M$ and $N$ are $p$-dimensional and $q$-dimensional \break
Euclidean manifolds respectively, i.e. $M\simeq {\bf    R}^p$ and
$N\simeq {\bf    R}^q$. The $total \; space$ will be the Euclidean
space $E=M\times N\simeq {\bf    R}^{p+q}$ coordinatized by the
independent and dependent variables $(x^1,...,x^p,u^1,...,u^q)$.
\end{Definition}
Let $f$ is a scalar-valued function with $p$ independent and $q$
dependent variables then it has $p_k={p+k-1\choose k}$ different
kth order partial derivatives and suppose the function
characterized by $u=f(x)$.
\begin{Definition}
Let $N^{(n)}$ is an appropriate space for representation of the
nth Tailor expansion of $f$, thus $N^{(n)}$ has the following
decomposition $$N^{(n)}\simeq {\bf    R}^q \times {\bf    R}^{p_1
q}\times {\bf    R}^{p_2 q}\times \cdots \times {\bf    R}^{p_n
q},$$ then $N^{(n)}$ is a $q{p+n\choose n}$ dimensional manifold.
Define
$$J^{(n)}=J^{(n)}E=M\times N^{(n)},$$ then $J^{(n)}$ is called the
$nth \; jet \; space$ of the total space $E$, which is a
$p+q{p+n\choose n}$ dimensional Euclidean manifold with vector
bundle structure whose fibers have $q{p+n\choose n}$ dimension.
\end{Definition}

According to the last definition, the nth jet space of ${\bf R}^2$
has dimension $n+2$ since we can write it in the form of ${\bf
R}^2={\bf    R}\times {\bf    R}$, ($M\simeq {\bf    R}$ and
$N\simeq {\bf    R}$).
\begin{Theorem}
Suppose G be a Lie group acting on E and {\fam=0\bf   v} is a
member of Lie algebra. Let $\Phi_x :G\rightarrow E$ be a smooth
map defined by $\Phi_x (g)=g.x$, then the infinitesimal generator
$\widetilde{v}$ corresponding v is given by $${\Phi_x}_* (v\mid_g
)=\widetilde{v}\mid_{g.x},$$ where ${\Phi_x}_*$ is the push
forward map of $\Phi_x$.
\end{Theorem}
See [2] for a proof.

Let $G$ be a Lie group acting on the total space $E$, then all of
infinitesimal generators of the group action have the following
form $$v=\sum_{i=1}^{p} \xi^{i} (x,u)\frac{\partial}{\partial
x^i}+\sum_{\alpha =1}^{q}
\varphi^{\alpha}(x,u)\frac{\partial}{\partial u^\alpha}\;.$$
\begin{Definition}
$Nth \; prolonged$ of $v$ denoted by $v^{(n)}$ is a vector field
on the $J^{(n)}E$, which is an infinitesimal generator of nth
prolonged action of $G$ on $E$.

$Nth \; prolonged$ action of $f$ represented by
$u^{(n)}=f^{(n)}(x)$, is a function from $M$ to $N^{(n)}$ (also
known as the $n-jet$ and denoted by $j_n f$), defined by
evaluating all the partial derivatives of $f$ up to order $n$.
\end{Definition}
Now it's time to give a formula to calculate nth prolongations.
The following theorem gives an explicit formula for the prolonged
vector field.

First of all we need two important definition.
\begin{Definition}
The $characteristic$ of a vector field $v$ on $E$ is a $q$-tuple
of functions $Q(x,u^{(1)})$, depending on $x$ and $u$ and first
order derivatives of $u$, defined by
$$Q^{\alpha}(x,u^{(1)})=\varphi^{\alpha}(x,u^{(1)})-\sum_{i=1}^{p}\xi^{i}(x,u)\frac{\partial u^{\alpha}}{\partial x^i}, \qquad \alpha=1,...q.$$
\end{Definition}
\begin{Definition}
Let $F(x,u^{(n)})$ be a differential function of order n. (A
smooth, real valued function $F:J^{(n)}\rightarrow {\bf    R}\;$,
defined on an open subset of the nth jet space is called a
$differential \; function$ of $order \; n$. ) The $total \;
derivative \; F$ with respect to $x^i$ is the (n+1)st order
differential function $D_i F$ satisfying
$$D_i [F(x,f^{(n+1)}(x))]=\frac{\partial}{\partial x^i}F(x,f^{(n)}(x)).$$
\end{Definition}
\begin{Theorem}
Let v be an infinitesimal generator on E, and let
$Q=(Q^1,...,Q^q)$ be its characteristic. The nth prolongation of
v is given explicitly by
$$v^{(n)}=\sum_{i=1}^{p} \xi^{i}(x,u)\frac{\partial}{\partial x^i}+\sum_{\alpha=1}^{q}\sum_{\sharp J=j=0}^{n} \varphi_{\alpha}^{J}(x,u^{(j)})\frac{\partial}{\partial u_{J}^{\alpha}}\;,$$
with coefficients
$$\varphi_{J}^{\alpha}=D_{J}Q^{\alpha}+\sum_{i=1}^{p}\xi^i u_{J,i}^{\alpha}\;,$$
where $J=(j_1,...,j_k)$ is a multi-indices where $1\leq j_k \leq
p$ and $\sharp J=k$.
\end{Theorem}
See [3; Theorem 2.36] for a proof.

\begin{Theorem}
(Infinitesimal Method, [2])A function $I:J^{(n)}\rightarrow {\bf
R}$ is a\break
 differential invariant for a connected transformation group
$G$ if and only if it is \break
annihilated by all prolonged
infinitesimal generators:
$$v^{(n)}(I(x,u^{(n)}))=0\;,\qquad v\in Lie(G)\;.$$
\end{Theorem}
Here $x=(x^1,...,x^p)$ is independent variable and $u^{(n)}$ is a coordinate chart on $N^{(n)}$.
\begin{Definition}
The collection $\{I_1,...,I_k\}$ of arbitrary differential
invariants is called $functionally\; independent$ if
$$dI_1\wedge ...\wedge dI_k \neq 0\;,$$ and $strictly\; functionally\; independent$
if $$dx\wedge du^{(n-1)}\wedge dI_1\wedge ...\wedge dI_k \neq
0\;.$$
\end{Definition}

Now we are ready to give some formulas for calculating $i_n$,
$j_n$, $s_n$ and $h_n$.

According to the prolonged action of $G$ on $J^{(n)}$ define
$V^{(n)}$ denoted a subset of $J^{(n)}$ which consists of all
points contained in the orbit of maximal dimension (generic
orbit). Then on $V^{(n)}$
$$i_n =\dim J^{(n)}-s_n =\dim J^{(n)}-\dim G+h_n\;,$$
thus the number of independent differential invariants of order
less or equal to n, forms a nondecreasing sequence $$i_0\leq
i_1\leq i_2\leq \cdots .$$

The difference $$j_n =i_n -i_{n-1},$$ is the number of strictly
functionally independent differential invariants of order n.

Note that $j_n$ cannot exceed the number of independent
derivative coordinate of order n, so if
$$q_n=\dim J^{(n)}-\dim J^{(n-1)}=q{p+n\choose n}-q{p-n-1\choose n-1}=q{p+n-1\choose n}\;,$$
is the number of derivative coordinate of order n, so $j_n\leq
q_n$, which implies that the elementary inequalities
$$i_{n-1}\leq i_n \leq i_{n-1}+q_n\;.$$

The maximal orbit dimension $s_n$ is also a nondecreasing
function of n, bounded by r, the dimension of $G$ itself:
$$s_0\leq s_1\leq s_2\leq \cdots \leq r\;.$$

On the other hand, since the orbit cannot increase in dimension
any more than the increase in dimension of the jet spaces
themselves, we have the elementary inequalities
$$s_{n-1}\leq s_n \leq s_{n-1} +q_n\;,$$ governing the orbit
dimension.
\section{ACTION OF LIE GROUP SL(2) AND SL(3)}
In this section we are going to define an action of each Lie group
on ${\bf    R}^2$, then we will give the infinitesimal generators
arise from Lie algebra and their prolongations, next we will
calculate $i_n, j_n, \cdots$.

For more details see [2], [3] and [6].
\subsection{Action of SL(2)}
Define an action of SL(2) on ${\bf    R}^2$,
$$(x,u)\longmapsto \Big(\frac{ax+b}{cx+d},u\Big)\; \mbox{where}\; \left(\begin{array}{cc}a&b\\c&d\end{array}\right)\in \mbox{SL}(2)\; \mbox{and}\; (x,u)\in {\bf    R}^2\;.$$
The Lie algebra of SL(2) generates by the following matrices:
\begin{eqnarray*}
A_1=\left(\begin{array}{cc}0&1\\0&0\end{array}\right)\qquad
A_2=\left(\begin{array}{cc}1&0\\0&-1\end{array}\right)\qquad
A_3=\left(\begin{array}{cc}0&0\\1&0\end{array}\right),
\end{eqnarray*}
and the infinitesimal generators corresponding to the above
matrices are:
\begin{eqnarray*}
{\fam=0\bf   v}_{1}=\frac{\partial}{\partial x}\;\;\;\;\;\;\;
{\fam=0\bf   v}_{2}=x\frac{\partial}{\partial x}\;\;\;\;\;\;\;
{\fam=0\bf   v}_{3}=x^2\frac{\partial}{\partial x}\;.
\end{eqnarray*}

According to the theorem 1.1.7, the three following vector fields
are the corresponding prolonged infinitesimal generators:
\begin{eqnarray*}
{\fam=0\bf   v}_{1}^{(n)} &=& \frac{\partial}{\partial x} \\
{\fam=0\bf   v}_{2}^{(n)} &=& x\frac{\partial}{\partial x}-\cdots -nu_{(n)}\frac{\partial}{\partial u_{(n)}} \\
{\fam=0\bf   v}_{3}^{(n)} &=& x^2\frac{\partial}{\partial
x}-\cdots
-\Big[n(n-1)u_{(n-1)}+2nxu_{(n)}\Big]\frac{\partial}{\partial
u_{(n)}}\;.
\end{eqnarray*}

According to the Theorem 1.1.8, the first three differential
invariants of order 0,1,2 are functions of $u$, thus
$i_0=i_1=i_2=1$, moreover the third one depends both on $u$ and
$$\frac{2u_{(1)}u_{(3)}-3u_{(1)}^2}{2u_{(1)}^4}\;,$$ thus $i_3=2$,
so the calculations show that $i_4=3,\cdots ,i_n=n-1$ and also
\break
 $s_0=1,s_1=2$ and $s_2=s_3=\cdots =s_n=3$, but $j_1=j_2=0$
and $j_0=j_3=\cdots =j_n=1$, so $h_0=2,h_1=1$ and $h_2=\cdots
=h_n=0$.
\subsection{Action of SL(3)}
Similarly, we define an action of SL(3) on ${\bf    R}^2$,
$$(x,u)\longmapsto \Big(\frac{ax+bu+c}{hx+ju+k},\frac{dx+eu+f}{hx+ju+k}\Big)$$
where
$\left(\begin{array}{ccc}a&b&c\\d&e&f\\h&j&k\end{array}\right)\in$
SL(3) and $(x,u)\in {\bf    R}^2$. The Lie algebra of SL(3)
generates by the following matrices:
\begin{eqnarray*}
A_1=\left(\begin{array}{ccc}0&0&1\\0&0&0\\0&0&0\end{array}\right)\qquad A_2=\left(\begin{array}{ccc}0&0&0\\0&0&1\\0&0&0\end{array}\right)
\end{eqnarray*}
\begin{eqnarray*}
A_3=\left(\begin{array}{ccc}1&0&0\\0&-1&0\\0&0&0\end{array}\right)\qquad A_4=\left(\begin{array}{ccc}0&0&0\\0&1&0\\0&0&-1\end{array}\right)
\end{eqnarray*}
\begin{eqnarray*}
A_5=\left(\begin{array}{ccc}0&0&0\\1&0&0\\0&0&0\end{array}\right)\qquad A_6=\left(\begin{array}{ccc}0&1&0\\0&0&0\\0&0&0\end{array}\right)
\end{eqnarray*}
\begin{eqnarray*}
A_7=\left(\begin{array}{ccc}0&0&0\\0&0&0\\1&0&0\end{array}\right)\qquad A_8=\left(\begin{array}{ccc}0&0&0\\0&0&0\\0&1&0\end{array}\right)\;,
\end{eqnarray*}
therefore the infinitesimal generators are:
\begin{eqnarray*}
\qquad \qquad \qquad \qquad {\fam=0\bf   v}_{1}  &=&
\frac{\partial}{\partial x} \\ {\fam=0\bf
v}_{2} &=& \frac{\partial}{\partial u} \\
\qquad \qquad \qquad \qquad {\fam=0\bf   v}_{3}  &=&
x\frac{\partial}{\partial x} \\ {\fam=0\bf
v}_{4} &=& u\frac{\partial}{\partial u} \\
\qquad \qquad \qquad \qquad {\fam=0\bf   v}_{5}  &=&
x\frac{\partial}{\partial u} \\ {\fam=0\bf
v}_{6} &=& u\frac{\partial}{\partial x} \\
\\ {\fam=0\bf   v}_{7}  &=& x^2\frac{\partial}{\partial
x}+xu\frac{\partial}{\partial u} \\ {\fam=0\bf   v}_{8} &=&
xu\frac{\partial}{\partial x}+u^2\frac{\partial}{\partial u}\;,
\end{eqnarray*}
and by calculating the prolonged infinitesimal generators, we
have,
\begin{eqnarray*}
{\fam=0\bf   v}_{1}^{(n)} &=& \frac{\partial}{\partial x} \\
{\fam=0\bf   v}_{2}^{(n)} &=& \frac{\partial}{\partial u} \\
{\fam=0\bf   v}_{3}^{(n)} &=& x\frac{\partial}{\partial
x}-\cdots -nu\frac{\partial}{\partial u_{(n)}} \\
{\fam=0\bf   v}_{4}^{(n)} &=& u\frac{\partial}{\partial
u}+\cdots +u_{(n)}\frac{\partial}{\partial u_{(n)}} \\
{\fam=0\bf   v}_{5}^{(n)} &=& x\frac{\partial}{\partial
u}+\frac{\partial}{\partial u_{(1)}} \\
{\fam=0\bf   v}_{6}^{(1)} &=& u\frac{\partial}{\partial
x}-u_{(1)}^2\frac{\partial}{\partial u_{(1)}} \\
{\fam=0\bf   v}_{6}^{(n)} &=& u\frac{\partial}{\partial x}-\cdots
-\Bigg\{\sum_{i=2}^{(n-1)/2}{n+1\choose i}u_{(i)}u_{(n-i+1)} \\
                           & & + \frac{1}{2}{n+1\choose
                           (n+1)/2}\Big[u_{((n+1)/2)}\Big]^2+(n+1)u_{(n)}\frac{\partial}{\partial
                           u_{(1)}}\Bigg\}\frac{\partial}{\partial
                           u_{(n)}}\;\; (n, \mbox{odd}) \\
{\fam=0\bf   v}_{6}^{(n)} &=& u\frac{\partial}{\partial x}-\cdots
-\Big[\sum_{i=2}^{n/2}{n+1\choose i}u_{(i)}u_{(n+i-1)} \\
                           & & +(n+1)u_{(n)}u_{(1)} \Big] \frac{\partial}{\partial
                           u_{(n)}}\;\;\;\;\;\; (n, \mbox{even}) \\
{\fam=0\bf   v}_{7}^{(n)} &=& x^2\frac{\partial}{\partial
x}+xu\frac{\partial}{\partial u}-\cdots
-\Big[n(n-2)+(2n-1)xu_{(n)}\Big]\frac{\partial}{\partial u_{(n)}} \\
{\fam=0\bf   v}_{8}^{(1)} &=& xu\frac{\partial}{\partial
x}+u^2\frac{\partial}{\partial
u}-(xu_{(1)}-u)u_{(1)}\frac{\partial}{\partial u_{(1)}} \\
{\fam=0\bf   v}_{8}^{(2)} &=& xu\frac{\partial}{\partial
x}+u^2\frac{\partial}{\partial x}-
(xu_{(1)}-u)u_{(1)}\frac{\partial}{\partial u_{(1)}}
-3xu_{(1)}u_{(2)}\frac{\partial}{\partial u_{(2)}} \\
{\fam=0\bf   v}_{8}^{(n)} &=& xu\frac{\partial}{\partial
x}+u^2\frac{\partial}{\partial u}-\cdots \\
                           & & - \Bigg\{{n\choose
                           (n+1)/2}\left[(n-2)u_{((n-1)/2)}+u_{((n+1)/2)}\right]u_{((n+1)/2)} \\
                           & & + \sum_{i=1+(n+1)/2}^{n}\left[{n\choose
                           j-1}(n-2)u_{(j-1)}+{n+1\choose j}
                               xu_{(j)}\right]u_{(i)}\Bigg\} \\
                           & & (j=1,...,\frac{n-1}{2}\;\;\mbox{and}\;\;n, \mbox{odd}) \\
{\fam=0\bf   v}_{8}^{(n)} &=& xu\frac{\partial}{\partial
x}+u^2\frac{\partial}{\partial u}-\cdots - \Bigg\{(\frac{n}{2}){n\choose n/2}[u_{(n/2)}]^2 \\
                           & & + \sum_{i=1+n/2}^{n}\left[{n+1\choose
                           j}xu_{(j)}+{n\choose
                           j-1}(n-2)u_{(j-1)}\right]
                           u_{(i)}\Bigg\}\frac{\partial}{\partial
                           u_{(n)}} \\
                           & & (j=1,...,\frac{n}{2}\;\;\mbox{and}\;\;n, \mbox{even})\;,
\end{eqnarray*}
corresponding to minimum value of $i$ in v$_8$$^{(n)}$, $j$ takes
its maximum value, i.e. when $i$ increases $j$ decreases.

According to the Theorem 1.1.8 with a tedious calculation we have,
$i_0=i_1=i_2=i_3=i_4=i_5=i_6=0$ and $i_7=1,i_8=2,\cdots ,i_n=n-6$,
consequently $j_0=j_1=j_2=j_3=j_4=j_5=j_6=0$ and $j_7=\cdots
j_n=1$, but $s_0=2,s_1=3,s_2=4,s_3=5,s_4=6,s_5=7$, thus
$h_0=6,h_1=5,h_2=4,h_3=3,h_4=2,h_5=1$,
 and $h_6=\cdots h_n=0$. The differential invariants have so
 complicated forms in this action. For example the first
 nonconstant differential invariant (the invariant of 7th prolonged
 action)is a function of
\begin{eqnarray*}
&& \hspace{-1cm}\frac{1}{18(40u_{(3)}^2+9u_{(2)}^2u_{(5)}-45u_{(2)}u_{(3)}u_{(4)})}\Big(11200u_{(3)}^8
\\
&& -33600u_{(2)}u_{(3)}^6u_{(4)}+6720u_{(2)}^2u_{(3)}^5u_{(5)} \\
&& +31500u_{(2)}^2u_{(3)}^4u_{(4)}^2-12600u_{(2)}^3u_{(3)}^3u_{(4)}u_{(5)}
\\
&& +720u_{(2)}^4u_{(3)}^3u_{(7)}-756u_{(2)}^4u_{(3)}^2u_{(5)}^2+\cdots
\\
&& -2835u_{(2)}^5u_{(4)}u_{(5)}^2-189u_{(2)}^6u_{(6)}^2\Big)\;,
\end{eqnarray*}
and the second differential invariant arise from 8th prolonged
action depends on the above phrase an the following one
\begin{eqnarray*}
&& \hspace{-1cm}\frac{1}{207360000}\Big((20412u_{(6)}^3+6561u_{(5)}^2u_{(8)}-26244u_{(5)}u_{(6)}u_{(7)})u_{(2)}^9
\\
&& +(((131220u_{(6)}u_{(7)}-65610u_{(5)}u_{(8)})u_{(4)}+104976u_{(5)}^2u_{(7)}
\\
&& \vdots \\
&& \mbox{ more than 10 lines calculation} \\
&& \vdots \\
&& -201600000u_{(2)}u_{(3)}^{10}u_{(4)}+4480000u_{(3)}^{12})/(u_{(3)}^3
\\
&& \frac{9}{40}u_{(2)}^2u_{(5)}-\frac{9}{8}u_{(2)}u_{(3)}u_{(4)})^4\Big)\;.
\end{eqnarray*}
\section{MORE DETAILS}
In this section we will give some details where in special cases we can classify the differential invariants.
\subsection{Differential Operators}
\begin{Theorem}
([2, Proposition 5.15])Suppose $E={\bf    R}\times N$ (N is a
q-dimensional manifold) be a total space. Let $I(x,u^{(n)})$ and
$J(x,u^{(n)})$ be functionally independent differential
invariants, at least one of which has order exactly n. Then the
derivative $\frac{dJ}{dI}=\frac{D_xJ}{D_xI}$ is an (n+1)st order
differential invariant.
\end{Theorem}
\begin{Definition}
If $I=I(x,u^{(n)})$ is any given differential invariant with one independent variable, then
$$\mathcal{D}=\frac{d}{dI}=(D_xI)^{-1}D_x\;,$$ is called an $invariant\; differential\; operator$
for the prolonged group action.
\end{Definition}

Since if $J$ is any differential invariant, so is $\mathcal{D}J$.
Therefore, we can iterate $\mathcal{D}$, producing a sequence
$$\mathcal{D}^k J=\frac{d^k J}{dI^k}\qquad k=0,1,...\;,$$
of higher order differential invariants. The last theorem
classifies differential operator on ${\bf    R}^2$.
\begin{Theorem}
([2, Proposition 5.16])Suppose G is a group of point
transformations acting on the jet space corresponding to $E={\bf
R}\times {\bf R}$. Then, for some $n\geq 0$, there are precisely
two functionally independent differential invariants $I$ and $J$
of order $n$ (or less). Furthermore for any $k\geq 0$, a complete
system of functionally independent differential invariants of
order n+k is provided by $I,J,\mathcal{D}J,...,\mathcal{D}^k J$,
where $\mathcal{D}$ is the associated invariant differential
operator.
\end{Theorem}

%GATHER{Xbib.bib}   % For Gather Purpose Only
%GATHER{Paper.bbl}  % For Gather Purpose Only
\bibliographystyle{amsplain}
\bibliography{xbib}
\end{document}